\documentclass[12pt]{amsart}
\usepackage{amssymb,amsmath,graphicx}
\usepackage{fancybox}
\usepackage{color,xcolor}            

 \newtheorem{theorem}{Theorem}[section]

 \theoremstyle{definition}
 
 \theoremstyle{remark}
 
 \newtheorem{algorithm}[theorem]{Algorithm}
 
 \numberwithin{equation}{section}

\allowdisplaybreaks 

\begin{document}
\title[Parallel Method for PBE]{A Parallel Method for Population Balance Equations Based
on the Method of Characteristics }

\author{Yu Li, Qun Lin, Hehu Xie}
\address{LSEC, Institute
of Computational Mathematics,
Chinese Academy of Sciences, Beijing 100190, China}
\email{liyu@lsec.cc.ac.cn }

\address{LSEC, Institute
of Computational Mathematics,
Chinese Academy of Sciences, Beijing 100190, China}
\email{linq@lsec.cc.ac.cn }

\address{LSEC, NCMIS, Institute
of Computational Mathematics,
Chinese Academy of Sciences, Beijing 100190, China}
\email{hhxie@lsec.cc.ac.cn }

\thanks{This work is supported in part
by the National Science Foundations of China (NSFC 11001259, 2011CB309703 and 2010DFR00700)
and Croucher Foundation of Hong Kong Baptist University, the national Center
 for Mathematics and Interdisciplinary Science, CAS, the President Foundation of AMSS-CAS.
 The second author gratefully acknowledges the
support from the MBF-Project SimParTurS
under the grant 03TOPAA1 and the Institute for Analysis and Computational Mathematics,
Otto-von-Guericke-University Magdeburg.}



\date{}

\begin{abstract}
In this paper, we present a parallel scheme to solve the population balance equations
based on the method of characteristics and the finite element discretization. The application of the
method of characteristics transform the higher dimensional population balance equation into a series of
lower dimensional convection-diffusion-reaction equations which can be solved in a parallel way.
Some numerical results are presented to show the accuracy and efficiency.
\end{abstract}

\maketitle

\section{Introduction}\label{Section_Introduction}
In this paper, we propose a parallel scheme to solve the population balance equation (PBE)
based on the application of the method of characteristics and the finite element method. The PBEs
aries from the model of the industrial crystallization process
(see, e.g., \cite{JohnRolandMitkovaSuudmacherTobiskaVoigt,Mersmann_1,Mersmann_2} and the reference cited therein).
 Recently, more and more
researchers are interested in the numerical methods for PBEs (c.f. \cite{AhmedMatthiesTobiska,Ganesan,GanesanTobiska,JohnRolandMitkovaSuudmacherTobiskaVoigt}).
In PBEs, besides the normal space and time variables, the distribution of entities also depends
on their own properties which are referred to as internal coordinates. It is a high dimensional system of
equations which is a big challenge from the computational point of view. In order to overcome this difficulty,
we use the method of characteristics (c.f. \cite{Chen,Evans}) to transfer the original problem to a series of lower-dimensional
convection-diffusion-reaction
 problems which are defined on the characteristics curves and the spatial directions. Based on the data structure for the method
 of characteristics, a parallel implementation can be applied to do the simulation process that can
 improve the computational efficiency.

So far, there exists the alternating direction (operator splitting) method for the PBE by decomposing the original problem
into two unsteady subproblems of smaller complexity (see, e.g., \cite{AhmedMatthiesTobiska,Ganesan,GanesanTobiska}).
In the two subproblems, the ordering
of the data for the solution needs to be different since they are discretized in different direction (c.f. \cite{AhmedMatthiesTobiska}).
It is not so suitable for the parallel implementation and prevents the further improvement of the
computation efficiency for the PBE.

In the present paper, we use the method of characteristics to transform the PBE into
a series of convection-diffusion-reaction equations on the characteristic curves in each time step.
Then the finite element method is applied to solve the series of convection-diffusion-reaction problems.
Furthermore, based on the data structure of the numerical scheme, a parallel scheme is constructed to solve
  the PBE based on the distributed memory. Some numerical results are provided to check the
  efficiency of this parallel method.

The following of the paper will go as follows: Section \ref{Section_Model} introduces the model problem under consideration and
defines some notation. In Section \ref{Section_Characteristics}, we describe method of characteristics for solving
 the PBE. The finite element discretization for the PBE
is described in Section \ref{Section_FEM}. Then Section \ref{Section_Parallel} gives the parallel implementation way
for the full discrete
form of the PBE. The numerical results are given in Section \ref{Section_Numerical_Results} to validate the efficiency of
the numerical method proposed in this paper. Some concluding remarks are given in the last section.

\section{Model problem}\label{Section_Model}
Let $\Omega_{\mathbf x}$ be a simply connected domain in $\mathcal{R}^d$ $(d=2\ {\rm or}\ 3)$ with Lipschitz continuous
 boundary $\partial \Omega_{\mathbf x}$,
$\Omega_{\ell}=[\ell_{\rm min}, \ell_{\rm max}]\subset \mathcal{R}$ and $T>0$. The state of the individual particle in
the PBE equation may consists of the external coordinate $\mathbf{x}$ ($\mathbf x=(x_1,\cdots,x_d)$), denoting its position in the
 physical space, and the internal coordinate $\ell$, representing the properties of particles, such as size, volume, temperature etc..
A PBE for a solid process such as crystallization with one internal coordinate can be described by the following
partial differential equation:

Find $z:(0,T]\times\Omega_{\ell}\times\Omega_{\mathbf x}\rightarrow\mathcal{R}$ such that
\begin{equation}\label{PBE}
\left\{
\begin{array}{ll}
\frac{\partial z}{\partial t}+G(\ell)\frac{\partial z}{\partial\ell}-\varepsilon\Delta_{\mathbf x} z
+\mathbf{b}(\mathbf{x})\cdot\nabla_{\mathbf x} z = f(t,\ell,\mathbf x) &
{\rm in}\ (0,T]\times\Omega_{\ell}\times\Omega_{\mathbf x}, \\
z(0,\ell,\mathbf x)=z_{\rm init}(\ell,\mathbf x) & {\rm in}\ \Omega_{\ell}\times\Omega_{\mathbf x}, \\
z(t,\ell_{\rm min},\mathbf x)=z_{\rm bdry}(t,\mathbf x) & {\rm on}\ (0,T]\times\Omega_{\mathbf x}, \\
z(t,\ell,\mathbf x)=0 & {\rm on} \ (0,T]\times \Omega_{\ell}\times \partial\Omega_{\mathbf x},
\end{array}
\right.
\end{equation}
where the diffusion coefficient $\varepsilon>0$ is a given constant, $\Delta_\mathbf{x}$ and $\nabla_\mathbf{x}$
denote the Laplacian and gradient with respect to $\mathbf{x}$, respectively, $\mathbf{b}$ is a given velocity and
satisfies $\nabla_\mathbf{x}\cdot\mathbf{b}=0$, and $f$ is a source function. Here $G(\ell)>0$ represents the growth
rate of the particles that depends on $\ell$ but is independent of $\mathbf{x}$ and $t$. Furthermore, let us assume
 the data  $G(\ell)$, $\mathbf b$, $f$, $z_{\rm init}$ and  $z_{\rm bdry}$ are sufficiently smooth functions for
 our error estimate analysis.

Now we introduce some notation of the function spaces (c.f. \cite{Chen,Ciarlet}). Let $H^m(\Omega_{\mathbf x})$
denote the standard Sobolev space of functions with derivatives up to $m$ in $L^2
(\Omega_{\mathbf x})$ and the norm is defined by
\begin{eqnarray*}
\|v\|_{H^m(\Omega_{\mathbf x})}&=&
\left(\int_{\Omega_{\mathbf x}}\sum^m_{0\leq |\alpha|\leq m}\Big|\frac{\partial^{\alpha} v}{\partial \mathbf x^{\alpha}}\Big|^2 d{\mathbf x}\right)^{1/2},
\end{eqnarray*}
where $\alpha$ denote a non-negative multi-index $\alpha=\{\alpha_1,\cdots,\alpha_d\}$, $|\alpha|=\sum_{1\leq j\leq d}\alpha_{j}$ and
\begin{eqnarray*}
\frac{\partial^{\alpha}v}{\partial^{\alpha}\mathbf x}&=& \frac{\partial^{\alpha_1\cdots\alpha_d}v}{\partial x_1^{\alpha_1}\cdots x_d^{\alpha_d}}.
\end{eqnarray*}
We use $(\cdot,\cdot)_{\mathbf x}$ and $\|\cdot\|_{L^2(\Omega_{\mathbf x})}$
to denote the $L^2$-inner product and the associated norm in $\Omega_{\mathbf x}$, respectively,
which are defined as follows
\begin{eqnarray*}
(v,w)_{\mathbf{x}}=\int_{\Omega_{\mathbf x}}vw d{\mathbf x}&{\rm and}&
\|v\|_{L^2(\Omega_{\mathbf x})}^2=(v,v)_{\mathbf x}.
\end{eqnarray*}
Let $X$ be a Banach space with the norm $\|\cdot\|_X$. Then we define
\begin{eqnarray*}
C(\Omega_{\ell};X)  &=&\Big\{v:\Omega_{\ell}\rightarrow X:\ v {\rm\ is\ continuous}\Big\}, \\
W^{m,\infty}(\Omega_{\ell};X)&=&\Big\{v: \Omega_{\ell}\rightarrow X:\ \Big\|\frac{\partial^j v}{\partial \ell^j}\Big\|_X<\infty,\
0\leq j\leq m\Big\},\\
W^{m,\infty}((0,T];X)&=&\Big\{v:(0,T]\rightarrow X:\ \Big\|\frac{\partial^jv}{\partial t^j}\Big\|_X <\infty,\ 0\leq j\leq m\Big\},
\end{eqnarray*}
where the derivatives $\partial^j v/\partial\ell^j$ and $\partial^j v/\partial t^j$ are understood in the sense of distribution on $\Omega_{\ell}$
and $(0,T]$, respectively.
The norms in the above defined spaces are given as follows
\begin{eqnarray*}
\|v\|_{C(\Omega_{\ell};X)}&=&\sup_{\ell\in\Omega_{\ell}}\|v(\ell)\|_X, \\
\|v\|_{W^{m,\infty}(\Omega_{\ell};X)}&=&\max_{0\leq j\leq m}\sup_{\ell\in\Omega_{\ell}}\Big\|\frac{\partial^jv}{\partial \ell^j}\Big\|_X,\\
\|v\|_{W^{m,\infty}((0,T];X)}&=&\max_{0\leq j\leq m}\sup_{t\in (0,T]}\Big\|\frac{\partial^jv}{\partial t^j}\Big\|_X.
\end{eqnarray*}
For spaces $X$, $Y$ and $Z$,  we use the short notation $Z(Y(X)):=Z((0,T];(Y(\Omega_{\ell};X))$ in this paper.

\section{Method of characteristics}\label{Section_Characteristics}
In this section, we describe the method of characteristics (c.f. \cite{Chen,Evans,Leveque})
for the PBE (\ref{PBE}). The reason we adopt this method for the discretization in the product space
$(0,T]\times \Omega_{\ell}$ is that it has the suitable data structure for the parallel implementation
which will be discussed in the following sections.

First we set
\begin{eqnarray*}
\psi(t,\ell)&=&(1+G(\ell)^2)^{1/2}.
\end{eqnarray*}
Let the characteristic direction associated with the hyperbolic part of (\ref{PBE}),
$\partial z/\partial t+G(\ell)\partial z/\partial \ell$, be denoted by $s(t)$. Then
\begin{eqnarray}
\frac{\partial }{\partial s} &=&\frac{1}{\psi}\frac{\partial}{\partial t}
+ \frac{G(\ell)}{\psi}\frac{\partial}{\partial \ell}.
\end{eqnarray}
Then (\ref{PBE}) can be written as
\begin{equation}\label{PBE_Characteristic}
\left\{
\begin{array}{ll}
\psi\frac{\partial z}{\partial s} -\varepsilon\Delta_{\mathbf x}z+\mathbf b(\mathbf x)\cdot\nabla_{\mathbf x}z=f&
{\rm in}\ (0,T]\times\Omega_{\ell}\times\Omega_{\mathbf x}, \\
z(0,\ell,\mathbf x)=z_{\rm init}(\ell,\mathbf x) & {\rm in}\ \Omega_{\ell}\times\Omega_{\mathbf x}, \\
z(t,\ell_{\rm min},\mathbf x)=z_{\rm bdry}(t,\mathbf x) & {\rm on}\ (0,T]\times\Omega_{\mathbf x}, \\
z(t,\ell,\mathbf x)=0 & {\rm on} \ (0,T]\times \Omega_{\ell}\times \partial\Omega_{\mathbf x},
\end{array}
\right.
\end{equation}
We use uniform partitions for the time interval $(0,T]$ and the internal coordinate interval $\Omega_{\ell}$, respectively.
Let $\tau=T/N$, $\iota=(\ell_{\rm max}-\ell_{\rm min})/M$, $t^n=n\tau$, $n=0,1,\cdots,N$ and
 $\ell_m=\ell_{\rm min}+m\iota$, $m=0,1,\cdots, M$.
In order to satisfy the stability condition, we set
\begin{eqnarray}\label{CFL_Condition}
\tau\leq\frac{\iota}{\max_{\ell_{\rm min}\leq \ell\leq \ell_{\rm max}}{G(\ell)}}.
\end{eqnarray}
Then starting with $z(0,\ell,\mathbf x)=z_{\rm init}$, $z(t,\ell_{\rm min},\mathbf x)=z_{\rm bdry}(t,\mathbf x)$,
 the equation (\ref{PBE_Characteristic}) can be discreted in each sub-intervals
 $(t^{n-1},t^n]\times(\ell_{m-1},\ell_m]\times\Omega_\mathbf{x}$ ($n=1,2,\cdots,N$ and $m=1,2,\cdots,M$) as follows.

First we compute
\begin{eqnarray}
\check{\ell}_m&=&\ell_m-\tau G(\ell_m).
\end{eqnarray}
Actually, this is a first order discretization to obtain the approximation at the time level $t=t^{n-1}$
for the following characteristic ordinary differential equation (c.f. \cite{Evans}):
\begin{equation}\label{ODE_Characteristic}
\left\{
\begin{array}{rcl}
\frac{d\ell}{dt}&=&G(\ell)\ \ \ \ {\rm in}\ [t^{n-1},t^n),\\
\ell(t^n)&=&\ell_m.
\end{array}
\right.
\end{equation}

From the condition (\ref{CFL_Condition}), we have $ \check{\ell}_m>=\ell_{\rm min}$ for $m\geq 1$.
Then we compute the direction differential $\psi\frac{\partial z}{\partial s}$ at the
node $(t^n,\ell_m)$ in the following way
\begin{eqnarray}
\psi(t^n,\ell_m)\frac{\partial z}{\partial s}(t^n,\ell_m,\mathbf x)&\approx&
\psi(t^n,\ell_m) \frac{z(t^n,\ell_m,\mathbf x)-\check{z}(t^{n-1},\check{\ell}_m,\mathbf x)}{(\tau^2+(\ell_m-\check{\ell}_m)^2)^{1/2}}\nonumber\\
&=&\frac{z(t^n,\ell_m,\mathbf x)-\check{z}(t^{n-1},\check{\ell}_m,\mathbf x)}{\tau},
\end{eqnarray}
where $\check{z}(t^{n-1},\check{\ell}_m,\mathbf x):=\alpha_m^nz(t^{n-1},\ell_{m-1},\mathbf x)+(1-\alpha_m^n)z(t^{n-1},\ell_{m},\mathbf x)$ with
$\alpha^n_m=\frac{\ell_m-\check{\ell}_m}{\iota}$.


In order to give the semi-discrete form of the PBE, we set $z_m^n(\mathbf x)\approx z(t^n,\ell_m,\mathbf x)$.
Then the semi-discrete form of the PBE can be defined as follows:
\begin{equation}\label{PBE_Semi_Discrete}
\left\{
\begin{array}{ll}
\frac{z_m^n(\mathbf x)-\check{z}_m^n(\mathbf x)}{\tau}-\varepsilon\Delta_{\mathbf x}z_m^n(\mathbf x)
+\mathbf b(\mathbf x)\nabla_{\mathbf x}z_m^n(\mathbf x)=f_m^n(\mathbf x)&{\rm in}\ \Omega_{\mathbf x},\\
z_m^0(\mathbf x)=z_{\rm init}^m(\mathbf x) &{\rm for}\ x\in \Omega_{\mathbf x}, \\
z_0^n(\mathbf x)=z_{\rm bdry}(t^n,\mathbf x) & {\rm for}\ (0,T]\times \Omega_{\mathbf x},\\
z_m^n(\mathbf x)=0 \ {\rm for}\ m=1,2,\cdots,M& {\rm on} \ \partial\Omega_{\mathbf x}
\end{array}
\right.
\end{equation}
where $f_m^n(\mathbf x)=f(t^n,\ell_m,\mathbf x)$,
$\check{z}_m^n(\mathbf x)=\alpha_m^nz_{m-1}^{n-1}(\mathbf x)+(1-\alpha_m^n)z_m^{n-1}(\mathbf x)$.

From the Taylor expansion method, we can derive the following error estimate for the semi-discrete form (\ref{PBE_Semi_Discrete})
\begin{eqnarray}\label{Error_Estimate_Semi_Discrete}
\|z(t^n,\ell_m,\mathbf x)-z_m^n(\mathbf x)\|_{C(C(X))}&\leq& C\tau\|z(t,\ell,\mathbf x)\|_{W^{2,\infty}(W^{1,\infty}(X))},
\end{eqnarray}
where the space $X$ can be $L^2(\Omega_{\mathbf x})$ or $H^1(\Omega_{\mathbf x})$.

\section{Finite element method}\label{Section_FEM}
In this section, we give the fully discrete form of the PBE by the finite element method.
Let $V_h$ be a finite element subspace of $H_0^1(\Omega_{\mathbf x})$ which has the $k$-th order accuracy (c.f. \cite{Chen,Ciarlet}):
\begin{eqnarray}
\inf_{v_h\in V_h}\|u-v_h\|_{H^1(\Omega_{\mathbf x})}
\leq Ch^{k}\|u\|_{H^{m+1}(\Omega_{\mathbf x})}\ \ \ \ \forall u\in H^{m+1}(\Omega_{\mathbf x}).
\end{eqnarray}
and
\begin{eqnarray}
\inf_{v_h\in V_h}\|u-v_h\|_{L^2(\Omega_{\mathbf x})}
\leq Ch^{k+1}\|u\|_{H^{m+1}(\Omega_{\mathbf x})}\ \ \ \ \forall u\in H^{m+1}(\Omega_{\mathbf x}).
\end{eqnarray}

Based on the finite element space $V_h$, we can define the fully discrete form for the PBE as follows:

For the $n$-th time step $t=t^n$ and $m=0, 1, \cdots, M$, find $z_{m,h}^n\in V_h$ such that
\begin{equation}\label{PBE_Full_Discrete_FEM}
\left\{
\begin{array}{ll}
\left(\frac{z_{m,h}^n-\check{z}_{m,h}^n}{\tau}, v_h\right)+
a(z_{m,h}^n,v_h)= (f_m^n(\mathbf x),v_h) & \forall v_h\in V_h,\\
a_0(z_{m,h}^0,v_h)= a_0(z_{\rm init}(\ell_m,\mathbf x),v_h) & \forall v_h\in V_h,\ m=1,\cdots,M,\\
a_0(z_{0,h}^n,v_h)= a_0(z_{\rm bdry}(t^n,\mathbf x),v_h) & \forall v_h\in V_h,
\end{array}
\right.
\end{equation}
where
$\check{z}_{m,h}^n =\alpha_m^nz_{m-1,h}^{n-1}+(1-\alpha_m^n)z_{m,h}^{n-1}$ with $\alpha_m^n$ being defined
in Section \ref{Section_Characteristics} and
\begin{eqnarray*}
a(u,v)&=&\int_{\Omega_{\mathbf x}}\big(\varepsilon\nabla u\cdot\nabla v + \mathbf b(\mathbf x)\cdot\nabla u \ v\big) d\mathbf x,\\
a_0(u,v)&=&\int_{\Omega_{\mathbf x}}\nabla u\cdot\nabla vd\mathbf x.
\end{eqnarray*}

From the standard error estimate theory of the finite element method (c.f. \cite{Chen,Ciarlet}),
 the fully discrete form (\ref{PBE_Full_Discrete_FEM}) has the following error estimates
\begin{eqnarray}\label{Error_Estimate_Full_Discrete_H1}
\max_{1\leq m\leq M}\|z(T,\ell_m,\mathbf x)-z_{m,h}^N\|_{H^1(\Omega_{\mathbf x})}\leq
C(\tau+h^k)\|z\|_{W^{2,\infty}(W^{1,\infty}(H^{k+1}(\Omega_{\mathbf x})))}
\end{eqnarray}
and
\begin{eqnarray}\label{Error_Estimate_Full_Discrete_L2}
\max_{1\leq m\leq M}\|z(T,\ell_m,\mathbf x)-z_{m,h}^N\|_{L^2(\Omega_{\mathbf x})} \leq
C(\tau+h^{k+1})\|z\|_{W^{2,\infty}(W^{1,\infty}(H^{k+1}(\Omega_{\mathbf x})))}.
\end{eqnarray}

\section{A parallel way}\label{Section_Parallel}

In this section, we present a parallel scheme to solve the PBE (\ref{PBE}) based on the
full discrete (\ref{PBE_Full_Discrete_FEM}). Fortunately, from (\ref{PBE_Full_Discrete_FEM}),
we can find the finite element equation is independent for each $m$ in any time step $t^n$.
Based on this property, we can construct a type of parallel scheme to implement the full
discretization of the fully discrete PBE (\ref{PBE_Full_Discrete_FEM}).

Assume we use $P$ processors to compute the PBE. Decompose the set $\{0,1,2,\cdots,M\}$ into $P$
 subsets $\mathbf m_1,\mathbf m_2,\cdots,\mathbf m_P$ such that
 $\mathbf m_1=\{m_0=0,1,\cdots,m_1-1\}$, $\mathbf m_p=\{m_{p-1},m_{p-1}+1,\cdots,
 m_{p}-1\}$ ($p=2, \cdots, P-1$) and $\mathbf m_P=\{m_{P-1},\cdots,m_P-1=M\}$.
 In the $p$-th processor, the equation (\ref{PBE_Full_Discrete_FEM}) is solved on the sub-intervals
 $(t^{n-1},t^n]\times (\ell_{m_{p-1}},\ell_{m_{p}-1}]\times\Omega_{\mathbf x}$
 ($n=1,2,\cdots,N$, $p=1,2,\cdots,P$, $\ell_{0}=\ell_{\rm min}$ and $\ell_{M}=\ell_{\rm max}$).
Because the growth rate of the particles $G(\ell)>0$, the dependence of each point $\ell_m$ is on the left ($\ell<\ell_m$),
the solution $z_{m_p-1,h}^{n-1}$ in the $p$-th processor as the initial condition for the $p+1$-th processor computing
at the $t^n$ time step.

We allocate the memory in the $p$-th processor ($p=1,\cdots,P$) to save the solutions $z_{m_{p-1},h}^n,\cdots,z_{m_{p}-1,h}^n$
and the $p$-th processor ($p=1,\cdots,P-1$) should send its saved solutions to the next $p+1$-th processor
after each time step computation.
Obviously, for $p=1$, we need to use the boundary condition $z_{\rm bdry}(t,\mathbf x)$ and
for $p=P$, the sending of solutions is not required since it is the last processor.
Based on this distribution of the memory and the computation of the scheme (\ref{PBE_Full_Discrete_FEM}),
we can construct the following parallel algorithm for the PBE.

%

\begin{algorithm}\label{Parallel_Scheme}
Parallel algorithm for PBE

For $n=1,2,\cdots,N$ Do
\begin{enumerate}
\item On each processor, compute the solution $z_{m,h}^n$ for $m\in \mathbf m_p$ ($p=1,2,\cdots,P$)
in this sub-interval $(t^{n-1},t^n]\times (\ell_{m_{p-1}},\ell_{m_p-1}]$.
\item For $p=1,2,\cdots,P-1$, send the solutions in the $p$-th processor $z_{m,h}^n\ (m\in\mathbf m_p)$
to the $p+1$-th processor.
\item If $n<N$, set $n:=n+1$ and go to Step 1. Else stop.
\end{enumerate}
\end{algorithm}

\section{Numerical results}\label{Section_Numerical_Results}
In this section, we provide some numerical results to validate the numerical scheme
proposed in this paper.
Let $\Omega_{\mathbf x} = [0,1]\times[0,1]$, $\Omega_{\ell}=[0,1]$, $T = 1$,
$\varepsilon=1$ and $\mathbf b(\mathbf x)= (1,1)^T$. We chose the functions
$f(t,\ell,\mathbf x)$, $z_{\rm init}(\ell,\mathbf x)$ and $z_{\rm bdry}(t,\mathbf x)$
 such that the exact solution is
\begin{eqnarray*}
z(t,\ell,x,y) = e^{-at}\sin(\pi \ell)\sin(\pi x)\sin(\pi y)
\end{eqnarray*}
with $a = 0.1$. The growth rate of the particles is $G(\ell)=\frac{1}{2}+2(1-\ell)\ell$.

First, we check the convergence order for the error estimates
\begin{eqnarray*}
\|e\|_0=\max_{1\leq m\leq M}\|z(T,\ell,\mathbf x)-z_{m,h}^n\|_{L^2(\Omega_{\mathbf x})}
\end{eqnarray*}
and
\begin{eqnarray*}
\|e\|_1=\max_{1\leq m\leq M}\|z(T,\ell,\mathbf x)-z_{m,h}^n\|_{H^1(\Omega_{\mathbf x})}.
\end{eqnarray*}

The convergence order of the linear  and quadratic finite element method for the discretization in $\Omega_{\mathbf x}$ is shown
in Tables \ref{Order_FEM_P1} and  \ref{Order_FEM_P2}.
From Tables \ref{Order_FEM_P1} and \ref{Order_FEM_P2}, we can find the numerical method with the linear and quadratic
finite element method in the space direction has the reasonable convergence order.

\begin{table}[htb!]
\caption{Error and rate of convergence in the space direction for $P1$
with $\tau=\iota=h^2$}\label{Order_FEM_P1}
\smallskip
\centering
\begin{tabular}{c|l|l}
\hline
mesh size $h$ & $\ \ \ \ \ \ \ \ \ \|e\|_0$  &\ \ \ \ \ \ \ \ \ $\|e\|_1$ \\
&error\ \ \ \ \ \ \ \ \ \ \ order &error\ \ \ \ \ \ \ \ \ \ \ order\\
\hline
$2^{-2}$  & 4.5702E-01  \ \           &2.6897E-00 \ \        \\
$2^{-3}$  & 1.4872E-01  \ \    1.6197 &1.5128E-00 \ \  0.8302 \\
$2^{-4}$  & 4.0481E-02  \ \    1.8773 &7.8083E-01 \ \  0.9541 \\
$2^{-5}$  & 1.0318E-02  \ \    1.9721 &3.9359E-01 \ \  0.9883 \\
$2^{-6}$  & 2.7230E-03  \ \    1.9219 &1.9720E-01 \ \  0.9970 \\
\hline
\end{tabular}
\end{table}
\begin{table}[htb!]
\caption{Error and rate of convergence in the space direction for $P2$
with $\tau=\iota=h^3$}\label{Order_FEM_P2}
\smallskip
\centering
\begin{tabular}{c|l|l}
\hline
mesh size $h$ & $\ \ \ \ \ \ \ \ \ \|e\|_0$  &\ \ \ \ \ \ \ \ \ $\|e\|_1$ \\
&error\ \ \ \ \ \ \ \ \ \ \ order &error\ \ \ \ \ \ \ \ \ \ \ order\\
\hline
$2^{-1}$  & 6.0137E-01  \ \           &2.5073E-00 \ \         \\
$2^{-2}$  & 6.3958E-02  \ \    3.2331 &8.5316E-01 \ \  1.5552 \\
$2^{-3}$  & 7.4660E-03  \ \    3.0987 &2.3528E-01 \ \  1.8584 \\
$2^{-4}$  & 9.5200E-04  \ \    2.9713 &6.0522E-02 \ \  1.9588 \\
\hline
\end{tabular}
\end{table}

We also check the convergence order for the method of characteristics developed in Section \ref{Section_Characteristics}.
The corresponding numerical result are provided in Table \ref{Order_Internal}. From this table, we can find the
convergence order is $1$ which is the same as in (\ref{Error_Estimate_Semi_Discrete}).
\begin{table}[htb!]
\caption{Error and rate of convergence in the internal coordinate for the method of
characteristics with $P2$ ($h=\iota$)} \label{Order_Internal}
\smallskip
\centering
\begin{tabular}{c|l|l}
\hline
mesh size $h$ & $\ \ \ \ \ \ \ \ \ \|e\|_0$  &\ \ \ \ \ \ \ \ \ $\|e\|_1$ \\
&error\ \ \ \ \ \ \ \ \ \ \ order &error\ \ \ \ \ \ \ \ \ \ \ order\\
\hline
$2^{-2}$ & 6.3862E-01 \ \          & 2.8423E-00 \ \           \\
$2^{-3}$ & 3.4562E-01 \ \   0.8858 & 1.5382E-00 \ \   0.8858  \\
$2^{-4}$ & 1.7650E-01 \ \   0.9695 & 7.8427E-01 \ \   0.9718  \\
$2^{-5}$ & 8.8689E-02 \ \   0.9928 & 3.9404E-01 \ \   0.9930  \\
$2^{-6}$ & 4.4398E-02 \ \   0.9983 & 1.9726E-01 \ \   0.9980  \\
\hline
\end{tabular}
\end{table}

Now we come to check the efficiency of the parallel scheme of Algorithm \ref{Parallel_Scheme}.  For this aim, we set the
discretization parameters $h=2^{-8}$, $\tau=\iota=1/512$ and the linear finite element method is adopted. The consuming
time (in seconds) are shown in Table \ref{Time_Parallel_Strong}. From Table \ref{Time_Parallel_Strong}, we can find the
parallel scheme Algorithm \ref{Parallel_Scheme} has good expansibility.
\begin{table}[htb!]
\caption{Strong parallel test with $P1$ ($h=2^{-8}$), $\tau=1/512$ and $\iota=1/512$}\label{Time_Parallel_Strong}
\centering
\begin{tabular}{l|c c c c c c c }
\hline
number of processors & 8 & 16 & 32 & 64 & 128 \\
\hline
time (in  seconds)   & 28103.01 & 13555.03 & 6832.26 & 3708.71 & 1840.43\\
rate of speed up     & 1.00 & 2.07 & 4.11 & 7.57 & 15.26\\
\hline
\end{tabular}
\end{table}

We also check the consuming time in each processor for different scale in each processor. For each test, we
run $8$ time steps ($N=8$).  Tables \ref{Time_Parallel_Weak_1} and \ref{Time_Parallel_Weak_2}
show the corresponding consuming time (in seconds) for the average time and maximum time, respectively,
for all the processors. These two tables also show that Algorithm \ref{Parallel_Scheme} has good
parallel property.
\begin{table}[htb!]
\caption{Weak parallel test with $P1$ element ($h=2^{-8}$): average time in seconds}\label{Time_Parallel_Weak_1}
\centering
\begin{tabular}{c|c c c c c c}
\hline
number in $\ell$    & 1  & 2  & 4  & 8  & 16  \\
\hline
8                & 9.30 & 15.30& 27.51& 55.28& 116.42  \\
16               & 9.91 & 15.44& 28.44& 59.44& 117.24  \\
32               & 9.85 & 16.98& 32.02& 60.93& 118.89  \\
64               & 10.01& 17.28& 32.66& 63.88& 121.96  \\
128              & 10.21& 17.98& 33.55& 64.27& 127.63  \\
\hline
\end{tabular}
\end{table}
\begin{table}[htb!]
\caption{Weak parallel test with $P1$ element ($h=2^{-8}$): maximum time in seconds}\label{Time_Parallel_Weak_2}
\centering
\begin{tabular}{c|c c c c c c}
\hline
number in $\ell$    & 1  & 2  & 4  & 8  & 16       \\
\hline
8                & 11.19& 16.10& 27.60& 60.52& 120.26  \\
16               & 11.26& 16.43& 31.54& 61.36& 120.83  \\
32               & 12.73& 18.50& 35.29& 68.18& 131.98  \\
64               & 11.20& 19.63& 36.39& 75.43& 133.55  \\
128              & 12.86& 20.28& 38.01& 73.63& 146.01  \\
\hline
\end{tabular}
\end{table}
\section{Concluding remarks}\label{Section_Remarks}
In this paper, we are concerned with the parallel numerical method for the PBEs
with one internal coordinate posed on the domain $(0,T]\times\Omega_{\ell}\times\Omega_\mathbf{x}$
with the dimension $1+1+d$.
The parallel scheme is based on the method of characteristics and the finite element discretization.
Some numerical results are also provided in Section \ref{Section_Numerical_Results}
 to demonstrate the efficiency of the proposed method.

Here, for the simplicity of the description of the numerical method, we assume the diffusion
coefficient $\varepsilon$ is large enough such that the diffusion is dominated.
For the convection dominated case (c.f. \cite{AhmedMatthiesTobiska,MatthiesSkrzypaczTobiska,RoosStynesTobiska}),
we will combine the method of characteristics and the
stabilized finite element methods (c.f. \cite{AhmedMatthiesTobiska,Chen,RoosStynesTobiska,MatthiesSkrzypaczTobiska})
and this is our future work.
Furthermore, the parallel method should also be applied to the simulation of the industrial
crystallization process (c.f. \cite{Mersmann_1,Mersmann_2}) and other similar models (c.f. \cite{JohnRolandMitkovaSuudmacherTobiskaVoigt}).

\section*{Acknowledgements}
This work is supported in part
by the National Science Foundations of China (NSFC 11001259, 2011CB309703 and 2010DFR00700)
and Croucher Foundation of Hong Kong Baptist University, the national Center
 for Mathematics and Interdisciplinary Science, CAS, the President Foundation of AMSS-CAS.
 The second author gratefully acknowledges the
support from the MBF-Project SimParTurS
under the grant 03TOPAA1 and the Institute for Analysis and Computational Mathematics,
Otto-von-Guericke-University Magdeburg.


\begin{thebibliography}{1}

\bibitem{AhmedMatthiesTobiska}
N. Ahmed, G. Matthies, L. Tobiska:
Finite element methods of an operator splitting applied to population balance equations.
J. Comput. Appl. Math., \textbf{236} (2011), 1604--1621.

\bibitem{Chen}
Z. Chen: \emph{Finite Element Methods and Their Applications}.
Springer, 2005.

\bibitem{Ciarlet}
P. Ciarlet:
\emph{The Finite Element Method for Elliptic Problem}.
North-Holland Amsterdam, 1978.



\bibitem{Evans}
L. Evanns:
\emph{Partial Differential Equations}.
American Mathematical Society, 1998.

\bibitem{Ganesan}
S. Ganesan, Population balance equations, Streamline-Upwind Petrov-Galerkin finite
element methods, operator-splitting method, backward Euler scheme, error analysis,
Preprint 1531, WIAS, Berlin, 2010.

\bibitem{GanesanTobiska}
S. Ganesan and L. Tobiska: Implementation of an operator splitting finite element
method for high-dimensional parabolic problems. Preprint 11-04, Fakult\"{a}t f\"{u}r Mathematik,
Otto-von-Guericke-Universit\"{a}t Magdeburg, 2010.

\bibitem{JohnRolandMitkovaSuudmacherTobiskaVoigt}
V. John, M. Roland, T. Mitkova, K. Sundmacher, L. Tobiska, and A. Voigt:
Simulations of population balance systems with one internal coordinate using finite
element methods. Chem. Eng. Sci., \textbf{64} (2009),  733--741.

\bibitem{Koch}
J. Koch:
Effiziente Behandlung von Integraloperatoren bei populationsdynamischen Modellen.
PhD thesis, Otto-von-Guericke-Universit\"{a}t Magdeburg, Fakult\"{a}t f\"{u}r Mathematik,
2005.

\bibitem{Leveque}
R. Leveque:
\emph{Finite Volume Methods for Hyperbolic Problems}.
Cambridge University Press, 2002.

\bibitem{MatthiesSkrzypaczTobiska}
G. Matthies, P. Skrzypacz, and L. Tobiska:
A unified convergence analysis for
local projection stabilisations applied to the Oseen problem.
M2AN Math. Model. Numer.
Anal., \textbf{41} (2007), 713--742.

\bibitem{Mersmann_1}
A. Mersmann:
Batch precipitation of barium carbonate.
Chem. Eng. Process., \textbf{38}
(1993), 6177--6184.

\bibitem{Mersmann_2}
A. Mersmann:
Crystallization and precipitation.
Chem. Eng. Process., \textbf{38} (1999), 345--353.

\bibitem{RoosStynesTobiska}
H.-G. Roos, M. Stynes, and L. Tobiska:
Robust numerical methods for singularly
perturbed differential equations. vol. \textbf{24} of Springer Series in Computational Mathematics,
Springer-Verlag, Berlin, second ed., 2008. Convection-diffusion-reaction and
fow problems.



\end{thebibliography}


\end{document}